\documentclass{svjour3}
\usepackage[english]{babel}
\usepackage{xfrac} %usecommand sfrac
\usepackage[utf8x]{inputenc}
\usepackage[T1]{fontenc}
\usepackage[a4paper,top=3cm,bottom=2cm,left=3cm,right=3cm,marginparwidth=1.75cm]{geometry}

\usepackage{amsmath}
\usepackage{graphicx}
\usepackage[colorinlistoftodos]{todonotes}
\usepackage[colorlinks=true, allcolors=blue]{hyperref}

\usepackage{mathptmx}
\usepackage{helvet}
\usepackage{courier}
\usepackage{type1cm}
\usepackage{makeidx}         % allows index generation
\usepackage{graphicx}        % standard LaTeX graphics tool
\usepackage{multicol}        % used for the two-column index
\usepackage[bottom]{footmisc}% places footnotes at page bottom
\usepackage{url}
\usepackage{amssymb,amsmath,mathrsfs,esint}
\usepackage{stmaryrd} %\varoast
\usepackage{phonetic} %\thorn
\usepackage{paralist} %a enlever
\usepackage{cancel}
\usepackage{textcomp}%to go with listings for quotation marks
\usepackage[normalem]{ulem} % allow strikethrough with \sout{}
\usepackage{listings}

\newcommand{\ML}{{\sc Matlab}}

\renewcommand{\paragraph}{\subparagraph}

\renewcommand{\vec}[1]{\ensuremath{\bm{\mathbf{#1}}}}
\newcommand{\mat}[1]{\ensuremath{\bm{\mathbf{#1}}}}
\usepackage{bm}

\newcommand{\D}{\ensuremath{\mat{D}}}

\newcommand{\bary}{\ensuremath{\beta}}

\newcommand\floor[1]{\lfloor#1\rfloor}

\newcommand{\DC}{\ensuremath{\mathbf{D}_{\mathrm{Chebyshev}}}}
\newcommand{\Dm}{\ensuremath{\mathbf{D}_{\mathrm{monomial}}}}
\newcommand{\DG}{\ensuremath{\mathbf{D}_{\mathrm{Degree-Graded}}}}
\newcommand{\DLeg}{\ensuremath{\mathbf{D}_{\mathrm{Legendre}}}}
\newcommand{\DL}{\ensuremath{\mathbf{D}_{\mathrm{Lagrange}}}}
\newcommand{\DN}{\ensuremath{\mathbf{D}_{\mathrm{Newton}}}}
\title{Differentiation Matrices for Univariate Polynomials}
\author{Amirhossein Amiraslani${}^1$, Robert M. Corless${}^2$, Madhusoodan Gunasingam${}^2$}
\date{\today}
\institute{${}^1$STEM Department, The University of Hawaii-Maui College, Kahului, Hawaii, \textsc{USA}\\
\& Faculty of Mathematics, K. N. Toosi University of Technology, Tehran, \textsc{Iran}\\
${}^2$The Ontario Research Center for Computer Algebra\\
\& The School of Mathematical and Statistical Sciences\\
The University of Western Ontario\\
London, Ontario, \textsc{Canada}}
\begin{document}
\maketitle

\begin{abstract}
We collect here elementary properties of differentiation matrices for univariate polynomials expressed in various bases, including orthogonal polynomial bases and non-degree-graded bases such as Bernstein bases and  Lagrange \& Hermite interpolational bases.
\end{abstract}
\section{Introduction}

The transformation of the (possibly infinite) vector of coefficients $\mathbf{a}={\{a_k\}}_{k\geq0}$ in the expansion

\begin{equation}
f(x)=\sum_{k\geq0} a_k\phi_k(x)
\end{equation}

\noindent to the vector of coefficients $\mathbf{b} = \{b_k\}_{k\geq0}$ in the expansion

\begin{equation}
f'(x) = \sum_{k\geq0} a_k\phi_k(x)
\end{equation}

\noindent is, of course, linear because the differentiation operator is linear\footnote{We are assuming that $f(x)$ is a differentiable function and that the set $\{\phi_k\}_{k\geq0}$, which we will now sometimes collect in a row vector $\mathbf{\phi}$, is a complete basis. However, the bulk of this paper will be about finite $\mathbf{\phi}$ representing polynomials of degree at most~$n$. Note that we represent $f'(x)$ and $f(x)$ in the same basis.}. Here the $\phi_k(x)$ are univariate polynomials. The matrix representation of the linear transformation from $\mathbf{a}$ to $\mathbf{b}$, denoted $\mathbf{b}=\mathbf{D}_{\phi}\mathbf{a}$, is 

\begin{equation}
\mathbf{D}_{\phi}=[d_{ji}]
\end{equation}

\noindent where $d_{ij} = \frac{\partial\phi_i'(x)}{\partial\phi_j(x)}$  (note the transposition); that is, the $d_{ij}$ are from

\begin{equation}
\phi_i'(x)=\sum_{j\geq0}d_{ij}\phi_j(x)\>.
\end{equation}
$\mathbf{D}_{\phi}$ is called a \textsl{differentiation matrix}.
In vector form, $f(x)=\boldsymbol{\phi}(x)\mathbf{a}$ (where $\boldsymbol{\phi}=[\phi_0(x),\phi_1(x),\dots]$) so

\begin{align}
f'(x)=&\boldsymbol\phi'(x)\mathbf{a} \nonumber \\
   =&\boldsymbol\phi(x)\mathbf{D}_{\phi}\mathbf{a} \nonumber\\
   =&\boldsymbol{\phi}(x)\mathbf{b}\>.   
\end{align}
Alternatively, we might work with $f'(x) = \mathbf{b}^T\boldsymbol{\phi}^T(x)$ and in that case use the transpose of $\D$, in $\mathbf{b}^T = \mathbf{a}^T\D^T$.

The most familiar differentiation matrix is of course that of the monomial basis $\phi_k(x)=x^k$. The $4\times 4$ differentiation matrix, for polynomials of degree at most 3, is in this basis,
\begin{equation}
      \Dm = \left[
	\begin{array}{cccc}
	0 & 1 & 0& 0 \\
	0& 0 & 2 & 0 \\
	0 & 0 & 0 & 3\\
    0 & 0 & 0 & 0\\
	\end{array}
	\right].
\end{equation}
This generalizes easily to the degree $n$ case. %$\boldsymbol{D}_{monomial}\mathbf{a}=[a_1,2a_2,\dots,na_n,0]^T% giving $\mathbf{D}_{monomial }\mathbf{a} = [a_1, 2a_2, 3a_3, .... ]^T$. 
This operation is so automatic that it's only rarely realized that it even has a matrix representation. If we are truncating to polynomials of degree at most~$n$ then 
the finite
 matrix $\mathbf{D}_{\textrm{monomial}}$ is
defined by:

\begin{equation}\label{DM}
      \mathbf{D}_{\textrm{monomial}} = \left[
	\begin{array}{cccccc}
	0 & 1 & 0&\cdots\\
     0   & 0 & 2 & 0 &\cdots \\
      & & 0  & 3 & 0 &\\
	&  &	 & \ddots & \ddots & \\
	&  &	 &     & &	n\\
	&&&&& 0\\ 
	\end{array}
	\right],
\end{equation}
an $n+1$ by $n+1$ matrix. 
	
Differentiation matrices in other bases, such as the Chebyshev basis, Lagrange interpolational basis, or the Bernstein basis, are also useful in practice and we will see several explicit examples.
	
\subsection{Reasons for studying Differentiation Matrices}

Differentiation matrices are used in spectral methods for the numerical solution of ordinary and partial differential equations, going back to their implicit use by Lanczos with the Chebyshev basis\footnote{Actually, Lanczos used the generalized inverse, $\DC^+$ which turns out to be bidiagonal for the Chebyshev basis; this is simpler for hand computation. We will see this later in this section.}. They can be used for quadrature, especially Filon or Levin quadrature, for highly oscillatory integrands. The first serious study seems to be~\cite{Don(1995)}. One of the present authors is working on this now~\cite{CorlessTrivedi:2018}. See also~\cite{Weideman},~\cite{olver2013fast}, and Chapter 11 of~\cite{corless2013graduate}.

In this paper we study differentiation matrices that occur when using various polynomial bases.  We confine ourselves to using one fixed basis $\{\phi_k\}_{k\geq0}$ for both $f(x)$ and $f'(x)$, but sometimes there are advantages to using different bases for $f'$ than for $f$: see~\cite{olver2013fast} for an example. The reasons for using the Chebyshev basis or the Lagrange basis include superior conditioning of expressions for functions in those bases, and sometimes superior convergence. The reasons for studying general properties and basis-independent properties, as this paper does, include the power of abstraction and the potential to apply to the results of other bases perhaps  more suited to the problem of your current interest. Another purpose is to see the relationships among the various bases. 

It helps exposition to have some example bases in mind, in order to make the general theory intelligible and interesting, so we describe the differentiation matrices for a few polynomial bases in the next section.

\subsection{Example Differentiation Matrices}

Before we give examples, we repeat the following general observation:
The columns of $\mathbf{D_{\phi}}$ are the coefficients of the derivatives $\phi_k'$ expressed in the $ \{\phi_k\}$.
\\
\noindent\textsl{Proof}: If $b=\phi_k(x)$ then $\mathbf{b}=\mathbf{e}_k$ and $\mathbf{D}\mathbf{b}=\mathbf{d}_k$ the $k$-th column of $\mathbf{D}$; but $\mathbf{b}'(x)=\sum_{j=0}^n c_j\phi_j(x)$ for some $c_j$, and $d_k=[c_0, c_1, ..., c_n]^T$, by definition.
\\
\noindent\textsl{Corollary}: if $\{\phi_k\}$ is degree-graded (\textsl{i.e.} $\deg\phi_k=k$), then $\mathbf{D}$ is strictly upper triangular. 
This is not true, of course, if $\{\phi_k\}$ is not degree-graded, \textsl{e.g.} $\phi_k$ is a Bernstein, Lagrange or Hermite interpolational basis.

\subsubsection{Chebyshev Polynomials}
One of the first kinds of differentiation matrices to be studied was for Chebyshev polynomials, \textsl{i.e}. $T_0(x)=1$, $T_1(x) = x$ and $T_{k+1} = 2xT_k(x)-T_{k-1}(x)$; alternatively, $T_k(x) = \cos(k\arccos(x))$ on $-1\leq x \leq1$. See for instance~\cite{olver2013fast} or (more briefly) Chapter 2 of~\cite{corless2013graduate}. For a thorough and modern introduction with application to the Chebfun software project see~\cite{berrut2004barycentric}. The derivative of  $T_k(x)$ is explicitly given in terms of $T_0, T_1,...,T_{k-1}$ as a sum, in~\cite{olver2013fast} and as a Maple program in~\cite{corless2013graduate}.

\begin{equation}
 \frac{dT_k(x)}{dx}=
    \begin{cases} 
      0 & k=0 \\
      k(\frac{1+(-1)^{k-1}}{2})T_0+2k\sum_{j=0}^{\floor{\frac{k-1}{2}}} T_{k-1-2j}(x)& k\geq 1 \>.\\
   \end{cases}
\end{equation}
Here the notation $\lfloor x\rfloor$ means the floor of $x$, the largest integer not greater than~$x$.
From this formula we may construct the infinite differentiation matrix $\DC$,
defined by
\begin{equation}\label{DC}
      \DC = \left[
	\begin{array}{ccccccccc}
	0 & 1 & 0 & 3 & 0 & 5 & 0 & 7 & \cdots\\
	 & 0 & 4 & 0 & 8 & 0 & 12 & 0 & \cdots\\
	 &   & 0  & 6 & 0 & 10 & 0 & 14 &\cdots\\
	 &&& 0 & 8 & 0& 12 & 0 &\cdots\\
	 &&&& 0 & 10 & 0 & 14&\cdots \\
     &&&&  & 0 & 12 & 0 &\cdots\\
     &&&&&& 0 & 14 & \ddots \\
     &&&&&&& 0 & \ddots \\
          &&&&&&& & \ddots
	\end{array}
	\right].
\end{equation}
As we see, the matrix is strictly upper triangular, just as the monomial basis matrix was; this is because the degree of $T_k'$ is $k-1$. Finite order differentiation matrices for Chebyshev polynomials are merely truncations of this. For a recent application of this matrix to the solution of pantograph equations, see~\cite{YANG2018132}.

\medskip
\goodbreak
\textbf{Remark}

Lanczos thought that this was cumbersome, and preferred the more compact \textsl{antiderivative} formulation (see~\cite{corless2013graduate} pp 125-126)
\begin{equation}
	\int T_k(x)dx=\frac{1}{2(k+1)} T_{k+1}(x)-\frac{1}{2(k-1)}T_{k-1}(x)+\frac{k\sin{k\pi/2}}{k^2-1}\>,
\end{equation}
(giving a correct value $\frac{(T_2(x)+T_0(x))}{4}$ in the limit as $k\rightarrow1$; also $\int T_0(x)dx=T_1(x)).$
This allows a more simple transformation from the \textsl{derivative}

\begin{equation}
f'(x)=\sum_{k\geq0} b_kT_k(x)
\end{equation}
to its \textsl{antiderivative}
\begin{equation}
f(x)=\sum_{k\geq0} a_kT_k(x)
\end{equation}

\noindent by what we will see is a generalized inverse of $\DC$:

The infinite 
tridiagonal matrix $\DC^+$,
derived from equation (7), is except for the first row
\begin{equation}
     \DC^{+} = \left[
	\begin{array}{cccccc}
	0 & 0 &   &  &  &\\
	1 &  0 & -1/2 &  & &\\
	  & 1/4 & 0 & -1/4 & &\\
	  &  & 1/6 & 0 & -1/6 &  \\
	 &&& 1/8 & 0 & \ddots\\
	 &&&& 1/10 &\ddots\\
	\end{array}
	\right]\>.
\end{equation}

This matrix is tridiagonal (with $0$ on main diagonal). Here the first row is $0$, meaning that an arbitrary constant can be added to the integral. We will see that the first row and the final column of truncations of this matrix will not matter for antiderivatives of degree $n-1$ polynomials.

\subsubsection{Legendre Polynomials}

The Legendre polynomials $\{P_n\}_n$ satisfy, $P_0(x)=1, P_1(x)=x$, 
\begin{equation}
\int_{-1}^1P_n(x)P_m(x)dx=\frac{2}{n+1}[n=m]
\end{equation}
This is a combinatorial notation for what elsewhere is termed the Kronecker Delta function. Here $[n=m]$ is 1 when $n=m$ and 0 otherwise. This is called Iverson's convention in~\cite{knuth1992two}. For a discussion of the merit of this notation, see~\cite{knuth1992two}. The Legendre polynomials satisfy the three term recursion relation
\begin{equation}
(n+1)P_{n+1}-(2n+1)P_n+nP_{n-1}=0
\end{equation}

By inspection, the differentiation matrix for polynomials $p(x)=\sum_{k\geq 0}c_kP_k(x)$ is, if $p'(x)=\sum_{k\geq 0}d_kP_k(x),$
\begin{gather}\label{DL}
 \begin{bmatrix} d_0 \\ 
      d_1 \\
      d_2 \\
      \vdots 
      \end{bmatrix}
 =
  \begin{bmatrix}
   0 & 1 & 0 & 1 & 0 & 1 & 0 & 1 \cdots \\
   0 & 0 & 3 & 0 & 3 & 0 & 3 & 0 \cdots \\
   0 & 0 & 0 & 5 & 0 & 5 & 0 & 5 \cdots \\
   0 & 0 & 0 & 0 & 7 & 0 & 7 & 0 \cdots \\
   & & & & \ddots & \ddots & \ddots & \ddots &\\
   \end{bmatrix}
  \begin{bmatrix}
   c_ 0\\
   c_1 \\ 
   c_2 \\
   \vdots
  \end{bmatrix}
\end{gather}
and 
\begin{gather}
 \begin{bmatrix} K \\ 
      c_1 \\
      c_2 \\
      c_3\\
      \vdots 
      \end{bmatrix}
 =
  \begin{bmatrix}
   0 & -\frac{1}{3}  \\
   1 & 0 & -\frac{1}{5} \\
    & \frac{1}{3} & 0 & -\frac{1}{7} \\
    &  & \frac{1}{5} & 0 & -\frac{1}{9} \\
    & & & \frac{1}{7 }& 0 & \ddots \\
    & & & & \frac{1}{9} & \ddots \\
        & & & &  & \ddots \\
   \end{bmatrix}
  \begin{bmatrix}
   d_ 0\\
   d_1 \\ 
   d_2 \\
   \vdots
  \end{bmatrix}.
\end{gather}
Like the matrix for the Chebyshev polynomials, the generalized inverse of $\DLeg$ is tridiagonal. The simplicity of these matrices recommend them.

\subsubsection{General Differentiation Matrix for Degree-Graded Polynomial Bases}
Real polynomials $\{\phi_n(x)\}_{n=0}^{\infty}$ with $\phi_n(x)$
of degree $n$ which are orthonormal on an interval of the real
line (with respect to some nonnegative weight function)
necessarily satisfy a three-term recurrence relation (see Chapter
10 of \cite{Davis(1963)}, for example). These relations can be written in the form 
\begin{equation} \label{eq.rec}
 x \phi_j(x)=\alpha_j \phi_{j+1}(x) +\beta_j \phi_j(x) +
\gamma_j\phi_{j-1}(x), \quad \quad j=0,1,\ldots, \end{equation} where the
$\alpha_j,\;\beta_j,\;\gamma_j$ are real, $\alpha_j\ne0$, $\phi_{-1}(x)\equiv 0$
$\phi_0(x)\equiv 1$.

Besides orthogonal polynomials, one can easily observe that the standard basis and Newton basis also satisfy (\ref{eq.rec}) with $\alpha_j= 1,\;\beta_j= 0,\;\gamma_j= 0$ and $\alpha_j= 1,\;\beta_j= z_j,\;\gamma_j= 0$, respectively where the $z_j$ are the nodes.

\textbf{Lemma}: $\DG$ has the following structure:

\begin{equation}\label{DG1}
      \DG = \left[
	\begin{array}{cccc}
               0&0&\cdots&0\\
              & & &  \vdots\\
              & \Huge{\mathbf{Q}}&& \\
              & & &0
	\end{array}
	\right]^{T},
\end{equation}

where 

\begin{equation}\label{DG}
\mathbf{Q}_{i, j}= \left\{\begin{array}{cl}
\frac{i}{\alpha_{i-1}}, & i=j\\
\frac{1}{\alpha_{i-1}}((\beta_{j-1}-\beta_{i-1})\mathbf{Q}_{i-1, j}+ \alpha_{j-2}\mathbf{Q}_{i-1, j-1}+ \gamma_j \mathbf{Q}_{i-1, j+1}- \gamma_{i-1}\mathbf{Q}_{i-2, j}). & i>j
\end{array}\right.\end{equation}
Any entry of $\mathbf{Q}$, with a negative or zero index is not considered in the above formula.\\

\noindent\textsl{Proof}: We provide the sketch of proof here. The proof itself is straightforward, but time-consuming. Taking the derivative of (\ref{eq.rec}) with respect to $x$, we have
\begin{equation} \label{deq.rec}
 x \phi'_j(x)+ \phi_j(x)=\alpha_j \phi'_{j+1}(x) +\beta_j \phi'_j(x) +\gamma_j\phi'_{j-1}(x), \quad \quad j=0,1,\ldots,   
\end{equation}

We let $j= 0$ in (\ref{deq.rec}) and simplify to get $$\phi'_1(x)= \frac{1}{\alpha_0}\phi_0(x).$$ We then let $j= 1$ in (\ref{deq.rec}) and simplify using (\ref{eq.rec}) with $j= 0$ and the result from the previous step to get $$\phi'_2(x)= \frac{\beta_0-\beta_1}{\alpha_0\alpha_1}\phi_0(x)+ \frac{2}{\alpha_1}\phi_1(x).$$ If we continue like this, and write the results in a matrix-vector form, the pattern stated in (\ref{DG}) will emerge.
\\

We can now find the matrices that we have for the monomial basis in (\ref{DM}), Chebyshev basis (\ref{DC}) and Legendre basis (\ref{DL}) directly from (\ref{DG1}) simply by plugging in the specific values for the $\alpha_j$, $\beta_j$, and $\gamma_j$ for each of them.

Another important degree-graded basis of this kind is the Newton basis. In the simplest case, let a polynomial $P(x)$ be specified by the
data $\{ \left(z_j,{P}_j\right) \}_{j=0}^{n}$ where the $z_j$'s
are distinct. The \textsl{Newton polynomials} are then defined by setting
$N_0(x)=1$ and, for $k=1,\cdots,n,$
$$  N_k(x)=\prod_{j=0}^{k-1}(x-z_j)\>.   $$
Then we may express \begin{equation} P(x)= \left[
	\begin{array}{ccccc}
              a_0&a_1& \cdots &a_{n-1}&a_n
	\end{array}
	\right]\left[
	\begin{array}{ccccc}
              N_0(x)\\N_1(x)\\\vdots\\N_{n-1}(x)\\N_n(x)
	\end{array}
	\right]\>.
\end{equation}

For $j=0,\cdots,n$, the $a_j$ can be found by divided differences as follows.

\begin{equation}
a_j=[P_0,P_1,\cdots,P_{j-1}],
\end{equation}
where we have $[P_j]=P_j$, and
\begin{equation}\label{divdif}
[P_{i},\cdots,P_{i+j}]=\frac{[P_{i+1},\cdots,P_{i+j}]-[P_i,\cdots,P_{i+j-1}]}{z_{i+j}-z_i}.
\end{equation}
A similar expression is possible even if the $z_j$ are not distinct, if we use \textsl{confluent} divided differences.  We return to this later, but note that the Newton polynomials are well-defined for $z_j$ that are not distinct.  Indeed, if they are all equal, say $z_j = a$, we recover Taylor polynomials $(z-a)^{j-1}$.

If in (\ref{eq.rec}), we let $\alpha_j=1$, $\beta_j=z_j$
and $\gamma_j=0$, it will become the Newton basis. For $n= 4$, $\DN$, as given by (\ref{DG1}), has the following form.

\begin{small}\begin{equation}\label{Nex}\hspace*{-0.5cm} \DN=\left[
	\begin{array}{ccccc} 0&0&0&0&0\\1&0&0&0&0 \\z_{0}-z_{1}&2&0&0&0\\ (z_{0}-z_{2})(z_{0}-z_{1})& -2z_{2
}+z_{1}+z_{0}&3&0&0\\ (z_{0}-z_{3})(z_{0}-z_{2})(z_{0}-z_{1})& (z_{1}-z_{3})( z_{1}-2z_{2}+z_{0}) +(z_{0}-z_{2})(z_{0}-z_{1}) &-3z_{3}+z_{2}+z_{1}+z_{{0}}&4&0
	\end{array}
	\right]^{T}\end{equation}\end{small}

\subsubsection{Lagrange Bases}

\bigbreak
Differentiation matrices for Lagrange bases are particularly useful. See \cite{corless2013graduate}
Chapter 2 for a detailed derivation. We give a summary here to establish notation. We suppose that function values                $\rho_k$  are given at distinct nodes $\tau_k$ (that is, $\tau_k=\tau_i \Leftrightarrow i=k,$ for $0\leq k\leq n$). Then the barycentric weights $\beta_k$ are found once and for all from the partial fraction expansion 
\begin{equation}
\frac{1}{w(\mathbf{z})}=\frac{1}{\prod_{k=0}^n (z-\tau_k)}=\sum_{k=0}^n \frac{\beta_k}{z-\tau_k},
\end{equation} 
giving

\begin{equation}
\beta_k=\prod_{\substack{j=0\\ j\neq k\\}}^n(\tau_k-\tau_j)^{-1}\>.
\end{equation}

These can be computed in a numerically stable fashion \cite{olver2013fast}, and once this has been done, the polynomial interpolant can be stably evaluated either by the first barycentric form 
\begin{equation}
\rho(z)=w(z)\sum_{k=0}^n \frac{\beta_k\rho_k}{z-\tau_k}
\end{equation}
or the second,
\begin{equation}
\rho(z)=\frac{\sum_{k=0}^n \frac{\beta_k\rho_k}{z-\tau_k}}{\sum_{k=0}^n \frac{\beta_k}{z-\tau_k}}\>.
\end{equation}
See \cite{berrut2004barycentric} for details. Here we are concerned with the differentiation matrix
\begin{equation}
\DL=[{d}_{ij}]
\end{equation}
(as derived in many places, but for instance see the aforementioned Chapter 11 of \cite{corless2013graduate}).\\ We have that 
\begin{equation}
d_{ij}=\frac{\beta_j}{\beta_i(\tau_i-\tau_j)} \text{ for } i\neq j 
\end{equation}
and
\begin{equation}
d_{ii}=-\sum_{j\neq i}d_{ij}
\end{equation}
Construction of this matrix is an $O(n^2)$ process, and evaluation of the vector of polynomial derivatives $b$ by
\begin{equation}
\mathbf{b}=\DL\rho\>.
\end{equation}
is also an $O(n^2)$ process. Once this has been done, then $\rho'(z)$ can be evaluated stably by re-using the previously computed barycentric weights:
\begin{equation}
\rho'(z)=w(z)\sum_{k=0}^n \frac{\beta_kb_k}{z-\tau_k}\>.
\end{equation}
If the derivative is to be evaluated \textsl{very} frequently, it may be cost-effective to modify the weights and throw away one node. This is usually not worth the bother.

\par\medskip\noindent
\textsl{Example}(taken from chapter 11 in \cite{corless2013graduate}) Note that if $\tau=[-1, -\frac{1}{3},\frac{1}{3}, 1]$ then the differentiation matrix is 
\begin{equation}
      \DL= \left[
	\begin{array}{cccc}
	-11 & \phantom{-}18 & -9 & \phantom{-}2\\
	 -2  & -3 & \phantom{-}6 & -1\\
	 \phantom{-}1   & -6  & \phantom{-}3 & \phantom{-}2\\
	 -2  &  \phantom{-}9 & -18 & \phantom{-}11\\
	\end{array}
	\right],
\end{equation}
so,
\begin{equation}
      \DL^+= \frac{1}{360}\left[
	\begin{array}{cccc}
	-81 & -147 & -123 & -9\\
	 -41 & -53 & -77 & -31\\
	 \phantom{-}31 & \phantom{-}77  & \phantom{-}53 & -41\\
	 \phantom{-}9  &  \phantom{-}123 & \phantom{-}147 & \phantom{-}81\\
	\end{array}
	\right]\>.
\end{equation}

If instead $\tau=[-1, -1/2, 1/2, 1]$, 
then it follows that
 
\begin{equation}
       \DL= \frac{1}{6}\left[
	\begin{array}{cccc}
	-19 & 24 & -8 & 3\\
	 -6 & 2 & 6 & -2\\
	 2 & -6  & -2 & 6\\
	 -3 & 8 & -24 & 19\\
	\end{array}
	\right],
\end{equation}
and
\begin{equation}
       \DL^+ = \frac{1}{720}\left[
	\begin{array}{cccc}
	-94 & -347 & -293 & 14\\
	 94 & -193 & -247 & -14\\
	 14 & 247  & 193 & -94\\
	 -14 & 293 & 347 & 94\\
	\end{array}
	\right].
\end{equation}
These matrices were displayed explicitly to demonstrate that, unlike the degree-graded case, the differentiation matrices are full, and their properties not very obvious\footnote{The row sums are zero, by design: the constant function has a constant vector representation, and its derivative should be (must be) zero. This is why $D_{ii}$ is the negative sum of all other entires.}.
If $\tau=[1, i, -1, -i]$,
\begin{equation}
 \DL = \frac{1}{2}\left[
	\begin{array}{cccc}
	 3 & -1+i & -1 & -1-i\\
	 -1+i & -3i & 1+i & i\\
	 1 & 1+i & -3 & 1-i\\
	 -1-i & -i & 1-i & 3i\\
	\end{array}
	\right],
\end{equation}
and
\begin{equation}
 \DL^+ = \frac{1}{24}\left[
	\begin{array}{cccc}
	 11 & 4-3i & 5 & 4+3i\\
	 -3+4i & 11i & 3+4i & 5i\\
	 -5 & -4-3i & -11 & -4+3i\\
	 -3-4i& -i5 & 3-4i & -11i\\
	\end{array}
	\right],
\end{equation}
which again has no obvious pattern (but see Theorem 11.3 of~\cite{corless2013graduate}: at least the singular values are simple).

\subsubsection{Hermite Interpolational Bases}
A Hermite interpolational basis is likely to be a bit less familiar to the reader than the Lagrange basis. They can be derived from Lagrange bases by letting two or more distinct nodes  ``flow together'' (from whence the word confluency comes).  Many methods to compute  Hermite interpolational basis representations of polynomials are known that fit consecutive function values and derivative values (i.e. $f(\tau_i)$, $f'(\tau_i)/1!$, $\ldots$, $f^{(s_i-1)}(\tau_i)/(s_i-1)!$ are consecutive scaled values of the derivatives of $f$ at a particular node $\tau_i$, which is said to have confluency $s_i$, a non-negative integer).

Many people use divided differences to express polynomials that fit confluent data, but this does not result in a Hermite interpolational basis (as pointed out in an earlier section, we would instead call that a Newton basis).  
We can solve the Hermite interpolation problem using a Newton basis, which is a degree-graded basis, and its differentiation matrix can be found through equation~\eqref{DG1}. 

Let's assume that at each node, $z_j$, we have the value and the derivatives of $P(x)$ up to the $s_j$-th order. The nodes at which the derivatives are given are treated as extra nodes. In fact we pretend that we have $s_j+1$ nodes, $z_j$, at which the value is $P_j$ and remember that $\sum_{i=0}^{k-1}s_i=n+1-k$. As such, the first $s_0+1$ nodes are $z_0$, the next $s_1+1$ nodes are $z_1$ and so on.

Using the divided differences technique, as given by equation~\eqref{divdif}, to find the $a_j$, whenever we get $[P_j,P_j,\cdots,P_j]$ where $P_j$ is repeated $m$ times, we have
\begin{equation}
[P_j,P_j,\cdots,P_j]=\frac{P^{(m-1)}_j}{(m-1)!},
\end{equation}
and all the values $P'_j$ to $P^{(s_j)}_j$ for $j=0,\cdots,k-1$ are given. For more details see e.g.~\cite{LJR(1983)}.

For this confluent Newton basis, like the simple Newton basis, $\alpha_j= 1$, $\beta_j= z_j$, and $\gamma_j= 0$, but some of the $\beta_j$ are repeated. Other than that, the differentiation matrix can be found for this basis in a manner identical to $\DN$. This approach was used in~\cite{AFS2017} to find the coefficients of the Birkhoff interpolation.

There are other advantages to solving the Hermite interpolation problem by using divided differences, for low degrees; the derivative is then almost directly available, for instance, and one does not really \textsl{need} a differentiation matrix.

But there are numerical stability disadvantages to the confluent Newton basis.  The main one is related to the relatively poor conditioning of the basis itself, for high-degree interpolants. [This does not matter much if the degree is low.] The next most important disadvantage is that the condition number of the polynomial expressed in this basis can be different if a different ordering of the nodes is used (it is usually true that the Leja ordering is good, but even so the condition number can be bad).  See~\cite{Corless(2004)} for numerical experiments that confirm this.

Another well-known solution to the Hermite interpolation problem involves constructing a basis that generalizes the Lagrange property, where each basis element is $1$ at one and only one node, and zero at all the others, which allows a direct sum to give the desired interpolant.
One possible such definition (there are many variations) for a Hermite interpolational basis is to define it as a set of polynomials $H_{i,j}(z)$ with the index $i$ corresponding to the node indices, so if the nodes are $\tau_i$ with $0 \le i \le n$ then again for $H_{i,j}(z)$ we would have $0 \le i \le n$.  The second index $j$ looks after the confluency at each node: $0 \le j \le s_i - 1$.  Importantly, one needs consecutive derivative data at each node (else one has a Birkhoff interpolation problem~\cite{AFS2017}\cite{butcher2011polynomial}).
Then we have the property (again written with the Iverson convention)
\begin{equation}
    {H_{i,j}^{(k)}(\tau_\ell)} = [i=\ell][j=k];
\end{equation}
that is, unless \textsl{both} the node indices are the same and the derivative indices are the same,
the given derivative of basis polynomial is zero at the given node; if both the node indices and the derivative indices \textsl{are} the same, then the (scaled) Hermite basis element takes the value $1$.  Using this definition, one can write the interpolant as a linear combination of this Hermite interpolational basis: $p(x) = \sum_{i=0}^n\sum_{j=0}^{s_i-1} \rho_{i,j}H_{i,j}(x)$.

But there is a better way, that uses a stable partial fraction decomposition to get a collection of generalized barycentric weights $\bary_{i,j}$ that can be used to write down an efficient barycentric formula for evaluation of the polynomial.  To be explicit, form the generalized node polynomial
\begin{equation}
    w(z) = \prod_{i=0}^n (z-\tau_i)^{s_i}\>,
\end{equation}
which is exactly what you would get from the Lagrange node polynomial on letting each group of $s_i \ge 1$ distinct nodes flow together.
Then the barycentric weights from the partial fraction decomposition
of $1/w(z)$ must now account for the confluency:
\begin{equation}
    \frac{1}{w(z)} = \sum_{i=0}^n \sum_{j=0}^{s_i-1} \frac{\beta_{i,j}}{(z-\tau_i)^{j+1}}\>. \label{eq:genparfrac}
\end{equation}
We will speak of the numerical computation of these $\beta_{i,j}$ shortly.  Once we have them, we may simply write down barycentric forms of the polynomial that solves the Hermite interpolational problem: the first form is
\begin{equation}
   p(z) = w(z) \sum_{i=0}^n \sum_{j=0}^{s_i-1} \sum_{k=0}^j \frac{\beta_{i,j}\rho_{i,k}}{(z-\tau_i)^{j+1-k}}\>.
\end{equation}
This form is simple to evaluate, and, provided the confluencies are not too large, numerically stable.
This form can be manipulated into a second barycentric form by replacing $w(z)$ with the reciprocal of its partial fraction expansion, equation~\ref{eq:genparfrac}.  The second form allows scaling of the generalized barycentric weights, which can prevent overflow. 
Incidentally, this allows us to give explicit expressions for the $H_{i,j}$ above:
\begin{equation}
    H_{i,j}(z) = \sum_{k=0}^{s_i-1} \beta_{i,j+k} w(z)(z-\tau_i)^{-k-1}\>.
\end{equation}
[Equivalent expressions are given in the occasional textbook, but not all works on interpolation do so; the formula seems to be rediscovered frequently.]

Given this apparatus, it makes sense to try to directly find the appropriate values of the derivatives at the nodes directly from the given function values and derivative values at the nodes; that is, by finding the differentiation matrix.
Rather than give the derivation (a complete one can be found in chapter 11 of~\cite{corless2013graduate}) we point to both Maple code and Matlab code that implements those formulas, at \url{http://www.nfillion.com/coderepository/Graduate_Introduction_to_Numerical_Methods/} in \texttt{BHIP.mpl} and \texttt{genbarywts.m}, respectively.  Evaluation in \ML\ can be done with the code \texttt{hermiteeval.m}.  We give an example below. 

If a polynomial is known at three points, say $[-1,0,1]$, and the values of $p$, $p'$, and $p''/2$ are known at $-1$, while the values of $p$, $p'$, $p''/2$, and $p'''/6$ are known at $0$, and the values of $p$ and $p'$ are known at $1$, then the differentiation matrix is found to be
%
%rho := [[r0, r1, (1/2)*r2], [s0, s1, (1/2)*s2, (1/6)*s3], [t0, t1]]
%tau := [-1, 0, 1]
%p, gam, DD := BHIP(rho, tau, z, 'Dmat' = true)
%latex(DD)
%

\begin{equation}
     \left[ \begin {array}{ccccccccc} 0&1&0&0&0&0&0&0&0
\\ \noalign{\medskip}0&0&2&0&0&0&0&0&0\\ \noalign{\medskip}-{\frac{201
}{2}}&-{\frac{177}{4}}&-15&96&-60&24&-12&9/2&-3/4\\ \noalign{\medskip}0
&0&0&0&1&0&0&0&0\\ \noalign{\medskip}0&0&0&0&0&2&0&0&0
\\ \noalign{\medskip}0&0&0&0&0&0&3&0&0\\ \noalign{\medskip}{\frac{83}{
4}}&6&1&-24&12&-12&4&{\frac{13}{4}}&-1/2\\ \noalign{\medskip}0&0&0&0&0
&0&0&0&1\\ \noalign{\medskip}35&11&2&0&48&0&16&-35&11\end {array}
 \right] \>.
\end{equation}
Applying this to the vector of values known at the nodes gives us the values of 
$[p'(-1)$, $p''(-1)$, $p'''(-1)/2$,
$p'(0)$, $p''(0)$, $p'''(0)/2$ , $p^{(iv)}(0)/6$, $p'(1)$, $p''(1)]^T$, which describe $p'(z)$ on these nodes in the same way that $p(z)$ was described. 

Notice that some rows are essentially trivial, and just move known values into their new places.  Notice that the nontrivial rows will, when multiplied by vectors representing constants (that is, $[c, 0, 0, c, 0, 0, 0, c, 0]^T$) give the zero vector.  The nontrivial rows are constructed by recurrence relations from the generalized barycentric weights $\bary_{i,j}$, which are themselves merely the coefficients in the partial fraction expansion of the node polynomial.

There is more than one way to compute the generalized barycentric weights $\beta_{i,j}$.  The fastest way that we know is the algorithm of~\cite{schneider1991hermite}, which internally uses a confluent Newton basis.  Unfortunately, because it does so, it inherits the poor numerical stability of that approach.  The codes referred to above use a direct local Laurent series expansion method instead, as outlined in~\cite{Henrici(1964)} for instance; this method is slower but much more stable. As discussed in~\cite{corless2013graduate}, however, it becomes less stable for higher confluency and cannot be perfectly backward stable even for $s_i \ge 3$. 
We will see an example in section~\ref{sec:HermiteExample}.

\subsubsection{Bernstein Polynomials}
The Bernstein differentiation matrix is a tridiagonal matrix. Its entries are as follows:
		\begin{equation}
			[\mathbf{D_B}]_{i,j} = 
			\begin{cases}
				2i-n & i=j \\
				-i  & j = i-1\\
				n-i & j=i+1
			\end{cases}\>.
		\end{equation}
Here the row and column indices $i$ and $j$ run from $0$ to $n$.
For polynomials of degree at most $n=4$ expressed in the Bernstein basis, the matrix is explicitly
\begin{equation}
     \left[ \begin {array}{ccccc} -4&4&0&0&0\\ \noalign{\medskip}-1&-2&3&0
&0\\ \noalign{\medskip}0&-2&0&2&0\\ \noalign{\medskip}0&0&-3&2&1
\\ \noalign{\medskip}0&0&0&-4&4\end {array} \right] \>.
\end{equation}
This is slightly different to the differentiation formulation seen in the Computer-Aided Geometric Design literature (e.g.~\cite{farin2014curves}), in that we \textsl{preserve the basis} to express the derivative in, even though that derivative is (nominally only) one degree too high.  Degrees of polynomials expressed in Bernstein bases can be \textsl{elevated}, however, and when they are too high, they can be \textsl{lowered} or \textsl{reduced}.  Indeed finding the \textsl{actual} degree of a polynomial expressed in a Bernstein basis can be, if there is noise in the coefficients, nontrivial.  Here we simply keep the basis that we use to express $p(x)$, namely
\begin{equation}
    p(x) = \sum_{i=0}^n c_i B_i^n(x)
\end{equation}
where 
\begin{equation}
    B_i^n(x) = {n \choose i} x^i (1-x)^{n-i}\>.
\end{equation}
By explicit computation, we find that the first column of the differentiation matrix (containing $-n$ in the zeroth row and $-1$ in the first row) correctly expresses the derivative of $B^n_0(x)$:
\begin{align}
    -nB^n_0(x) - B^n_1(x) &= -n(1-x)^{n} - n x(1-x)^{n-1} \nonumber\\
    &= -n(1-x)^{n-1}(1-x + x) \nonumber \\
    &= \frac{d}{dx}B^n_0(x)\>.
\end{align}
Similarly, for $1 \le i \le n-1$,
\begin{align}
    (n-i+1)B^n_{i-1} + (2i-n)B^n_i - (i+1)B^n_{i+1} &= x^{i-1}(1-x)^{n-i-1}{n \choose i}\left( i(1-x)^2 + (2i-n)x(1-x) - (n-i)x^2 \right) \nonumber \\
    &= x^{i-1}(1-x)^{n-i-1}{n \choose i}\left( i-nx \right) \nonumber \\
    &= \frac{d}{dx}B^n_i(x)\>.
\end{align}
By the reflection symmetry of $B^n_n(x)$ with $B^n_0(x)$, the final column is also correct. 
\begin{remark}
As with the Lagrange polynomial bases, the pseudo-inverse of the Bernstein basis differentiation matrix is full. Also as with the Lagrange case, because $1 = \sum B^n_k(x)$ (that is, the Bernstein basis forms a partition of unity), application of the Bernstein differentiation matrix to a constant vector must return the zero vector and hence the row sums must be zero.
\end{remark}

\section{Basic Properties}

\textsl{Definition}: Let $\mathbf{X}_{\phi}^k$ be the vector of coefficients of $x^k$ in the basis $\phi$. That is, if
\begin{equation}
x^k =b_{k,0}\phi_0(x)+b_{k,1}\phi_1(x)+\dots+b_{k,n}\phi_n(x)\>, 
\end{equation}
then
\begin{equation}
\mathbf{X}_{\phi}^k=[b_{k,0},b_{k,1},\dots,b_{k,n}]^T\>.\
\end{equation}
Set $\mathbf{1}_{\phi}=\mathbf{X}_{\phi}^0$. Let $\mathbf{V}$ be the matrix whose $k$-th column (numbering from zero) is $\frac{1}{k!}\mathbf{X}_{\phi}^k$. 
\bigbreak

\subsection{Eigendecomposition of Differentiation Matrices}
Let $\mathbf{D}_{\phi}$ be the differentiation matrix for polynomials of degree at most $n$, expressed in the basis $\{\phi_k\}_{k=0}^n$. Note that if 

\begin{equation}
b=\rho_0\phi_0+\rho_1\phi_1+...+\rho_n\phi_n
\end{equation}
then
\begin{equation}
\rho'=b_0\phi_0+b_1\phi_1+...+b_n\phi_n\>.
\end{equation}
same basis; for degree graded basis, $b_n=0$. Then $\mathbf{D}_{\phi}\mathbf{p}=\mathbf{b}$ where 
\begin{align}
\mathbf{\rho}=\begin{bmatrix}
           \rho_{0} \\
           \rho_{1} \\
           \vdots \\
           \rho_{n}\>
         \end{bmatrix}
\end{align} and 

\begin{align}
\mathbf{b}=\begin{bmatrix}
           b_{0} \\
           b_{1} \\
           \vdots \\
           b_{n}
         \end{bmatrix}.
\end{align}
\textbf{Lemma}: $\mathbf{D}$ is nilpotent.\\

\noindent\textsl{Proof}: $\mathbf{D}^{n+1}\mathbf{\rho}(x)=0$ for every  polynomial of degree at most $n$; hence $\mathbf{D}^{n+1}\mathbf{p}(x)=0$ as required.\\

\noindent$\mathbf{Remark}$: Therefore all eigenvalues are zero.\\

\noindent\textsl{Proposition}
$\mathbf{D}_{\phi}\mathbf{V}=\mathbf{VJ}$, where 
\begin{equation}
      \mathbf{J}= 
	 \left[
	\begin{array}{cccc}
	0 & 1 & &\\
	 0  & 0 & 1  & \\
	  \vdots&\vdots&\ddots&1\\
	 0&&&0\\
	\end{array}
	\right],
\end{equation}
is the Jordan Canonical Form of the differentiation matrix $\mathbf{D}_{\phi}$.
\bigbreak
\noindent\textsl{Proof}: $\mathbf{D}_{\phi}(\frac{1}{k!}\mathbf{X}_{\phi}^k)=\frac{1}{(k-1)!}\mathbf{X}_{\phi}^{k-1}$ for $k\geq 1$, by construction. Moreover the columns $\mathbf{X}_{\phi}^k$ are linearly independent because the monomials $1,x,x^2,\dots,x^n$ and $\phi$ are a polynomial basis. Thus $\mathbf{V}$ is invertible.\\

\noindent\textbf{Remark} The isomorphism of the polynomial representation by coefficient vectors (of the basis $\phi$) is complete for addition, subtraction, differentiation, and scalar multiplication; but the representation of $p\cdot q$ is possible only if $\deg p+\deg q\leq n$. The multiplication rules are interesting as well; we get the usual Cauchy convolution for the monomial basis.

\bigbreak
\subsection{Pseudoinverse}
\textsl{Observation}   

As long as $\deg p<n$, anti-differentiation works by using the pseudo inverse; one then adds a constant times $\mathbf{1}_{\phi}$. Call this anti-differentiation matrix $\mathbf{S}$. Then we want $\mathbf{S1}_{\phi}=\mathbf{X}_{\phi},$ and $\mathbf{S}\frac{\mathbf{X}^{k-1}}{(k-1)!}=\frac{\mathbf{X}^k}{k!}$.\\
\bigbreak
\noindent Therefore, \begin{equation}\mathbf{S}\mathbf{V}(1:n)=[\mathbf{X},\frac{\mathbf{X}^2}{2},\dots,\frac{\mathbf{X}^{n-1}}{(n-1)!},\frac{\mathbf{X}^n}{n!}]\>,\end{equation}

%What happens with $S\frac{\mathbf{X}^n}{n!}$?
\noindent and thus,
\begin{equation}
\mathbf{VJ^+V^{-1}[V_n]=VJ^+[0, 0, ..., 0, 1]^T=[0,0, ...., 0]^T\>.}
\end{equation}
\textbf{Lemma}: \textsl{The Moore-Penrose pseudo-inverse of}
\begin{equation}
      \mathbf{J}=  \left[
	\begin{array}{cccc}
	0 & 1 & &\\
	 0  & 0 & 1  & \\
	  \vdots&\vdots&\ddots&1\\
	 0&&&0\\
	\end{array}
	\right],
\end{equation}
is 
\begin{equation}
      \mathbf{J}^+=\mathbf{J^T}= \left[
	\begin{array}{cccccc}
	0 & 0 & & & &\\
	 1  & 0 &  & & & \\
	  &1 & \ddots & & & \\
	 && \ddots && &\\
	 &&&&1&0\\
	\end{array}
	\right].
\end{equation}

\textsl{Proof}: We need to verify that $\mathbf{JJ^TJ=J},\mathbf{J^TJ=J}$ and that both $\mathbf{J^TJ}$ and $\mathbf{JJ^T}$ are symmetric. The last two are trivial. Computation shows 

\begin{equation}
      \mathbf{JJ^T=J^TJ}= \left[
	\begin{array}{ccccc}
	0 & 0 & 0  &&\\
	 0  & 1 & 0  & & \\
	    &0 &1&&\\
	    &&&\ddots&\\
	 0&&&&1\\
	\end{array}
	\right],
\end{equation}
 so 
\begin{equation}
      \mathbf{JJ^TJ}= \left[
	\begin{array}{cccc}
	0 & 1 & &\\
	 0  & 0 & 1  & \\
	  \vdots&\vdots&\ddots&1\\
	 0&&&0\\
	\end{array}
	\right] = \mathbf{J}\>.
\end{equation}
 Similarly $\mathbf{J^TJJ^T=J}$.\\
 
$\mathbf{Proposition}$:
The matrix $\mathbf{D}^+=\mathbf{VJ^+V^{-1}}$ is a generalized inverse of $\mathbf{D}$. 

\textsl{Proof}: It suffices to verify only the first two of the Moore-Penrose conditions: $\mathbf{D^+DD^+}=\mathbf{D^+}$ and $\mathbf{DD^+D}=\mathbf{D}$.

These follow immediately. Interestingly $\mathbf{D^+}$ is not (in general) a Moore-Penrose inverse: neither $\mathbf{D^+D}$ nor $\mathbf{DD^+}$ need be Hermitian. 

The matrix $\mathbf{V}$ in a Lagrange basis is 

\begin{equation}
      \mathbf{V} = \left[
	\begin{array}{ccccc}
	1 & \tau_0  & \cdots & \cdots & \frac{\tau_0^n}{n!} \\
	 1 & \tau_1&  & & \\
	  \vdots & \vdots  &&&\\
	 1& \tau_n &&& \frac{\tau_n^k}{n!}
	\end{array}
    \right]\>.
\end{equation}
This is the product of a Vandermonde matrix and 
\begin{equation}
\left[
	\begin{array}{ccccc}
	1 &  &  &&\\
	   & 1 &  & & \\
	    & & \frac{1}{2} &&\\
	 &&& \ddots &\\
	 &&&& \frac{1}{n!}
	\end{array}
	\right].
\end{equation}
 
This is likely to be extraordinarily ill-conditioned. However, this gives an explicit JCF for differentiation matrices on Lagrange bases.

\section{Accuracy and Numerical Stability}
There are several questions regarding numerical stability (and, unfortunately, the answers vary with the basis used, and with the degree).  For the orthogonal polynomial bases and the Bernstein bases, the differentiation matrices have integer or rational entries, and there are no numerical difficulties in constructing them, only (perhaps) with their use.  For the Lagrange and the Hermite interpolational bases, the (generalized) barycentric weights need to be constructed from the nodes, and then the entries of the differentiation matrix constructed from the weights.  In floating-point arithmetic, this can be problematic for some sets of nodes (especially equally-spaced nodes); higher-precision construction of the weights, or use of symmetries as with Chebyshev nodes, may be needed.  High or variable confluency can also be a difficulty.
Use of higher precision in construction of the barycentric weights and of the differentiation matrix may be worth it, if the matrix is to be used frequently.

For all differentiation matrices, there is the question of accuracy of computation of the polynomial derivative by matrix multiplication. In general, differentiation is infinitely ill-conditioned: the derivative of $f(x) + \varepsilon v(x)$ can be arbitrarily different to the derivative of $f(x)$.  However, if both $f$ and the perturbation are restricted to be \textsl{polynomial}, then the ill-conditioning is finite, and the absolute condition number is bounded by the norm of the differentiation matrix $\mathbf{D}$.  This is Theorem 11.2 of~\cite{corless2013graduate}, which we state formally below.
\begin{theorem}
If $f(x)$ and $\Delta f(x)$ are both polynomials of degree at most $n$, and are both expressed in a polynomial basis $\phi$, then
\begin{equation}
    \| \Delta f'(x) \| \le \| \mathbf{D}_{\phi} \| \| \Delta f (x) \| 
\end{equation}
where the norms  $\| \Delta f' \|$ and $\| \Delta f \|$ are vector norms of their coefficients in $\phi$ and the norm of the differentiation matrix is the corresponding subordinate matrix norm. 
\end{theorem}
One should check the norm $\|\D\|$ whenever one uses a differentiation matrix.
We remark that the norms of powers of $\D$ can grow very large. For instance, for the Bernstein basis of dimension $n+1$ we find\footnote{We have no proof, only experimental evidence; it should be possible to prove this but we have not done so.} that $\| \D^n \|_\infty = 2^n n!$.  The next power gives the zero matrix, of course. To give a sense of scale, we have $\| \D \|_\infty = 2n$ and hence this norm to the $n$th power is much larger yet, being $(2n)^n$ so a factor $n^n/n! \approx \exp(n)/\sqrt{2\pi n}$ larger.
As a corollary, from the results discussed in~\cite{Embree:HLA:2013} the $\varepsilon$-pseudospectral radius of the $n+1$-dimensional Bernstein $\D$ matrix must then at least be $(2^n n!)^{1/(n+1)} \varepsilon^{1/(n+1)} \sim 2n \varepsilon^{1/(n+1)}/e$ as $n \to \infty$, for any $\varepsilon > 0$. This implies that for large enough dimension, matrices very near to $\D$ will have eigenvalues larger than $1$ in magnitude.  We believe that similar results hold for other bases, indicating that higher-order derivatives are hard to compute accurately by using repeated application of multiplication by differentiation matrices (as is to be expected).

\subsection{A Hermite interpolational example \label{sec:HermiteExample}}
Consider interpolating the simple polynomial that is identically $1$ on the interval $-1 \le z \le 1$, using nodes with confluency three.  That is, at each node we supply the value of the function ($1$), the value of the first derivative ($0$), and the value of the second derivative divided by $2$, which is also in this case just $0$.  We consider taking $n+1$ nodes $\tau_j$ for $0 \le j \le n$, which gives us $3(n+1)$ pieces of data and thus a polynomial of degree at most $3n+2$.  We then plot the error $p(z)-1$ on this interval.  We also compute the differentiation matrix $\D$ on these nodes with this confluency, and multiply $\D$ by the vector containing the data for the constant function $1$.  This should give us an identically $0$ vector (call it $\vec{Z}$), but will not, because of rounding error.  We compute the infinity norm of $\vec{Z}$ and the infinity norm of the matrix $\D$.

\begin{figure}[h]
    \centering
    \includegraphics[width=0.5\textwidth]{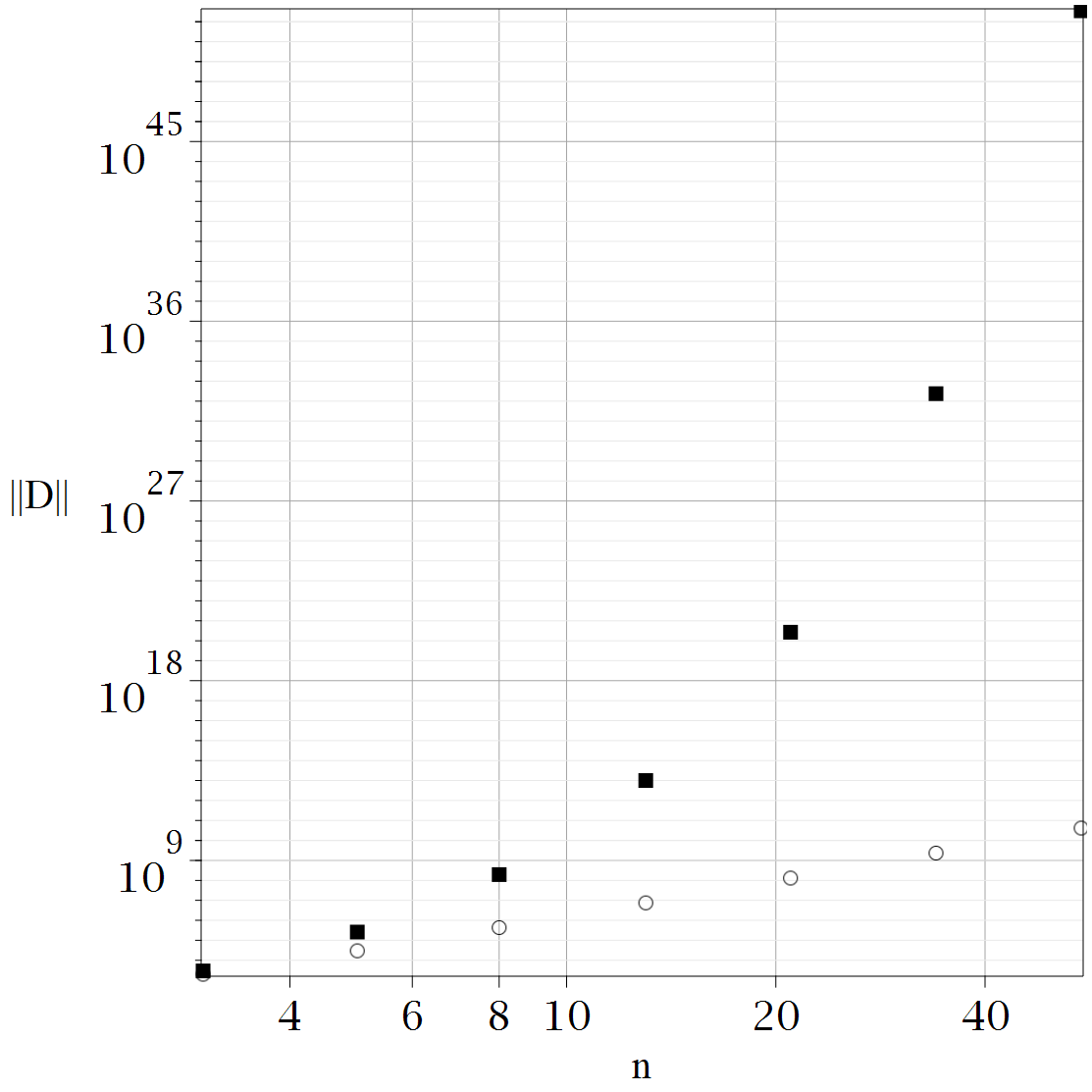}
    \caption{A comparison of norms of the differentiation matrices for Hermite interpolational basis on $n+1$ nodes, of confluency $3$, between equally-spaced nodes (solid boxes) and Chebyshev nodes (circles).  We see growth in $n$ for both sets of nodes, but much more rapid growth for equally-spaced nodes.}
    \label{fig:Dnorm}
\end{figure}

We take two sets of nodes: first the Chebyshev nodes $\tau_j = \cos(\pi(n-j)/n)$, and second the equally-spaced nodes $\tau_j = -1 + 2j/n$. We take $n=3$, $5$, $8$, $\ldots$, $55$ (Fibonacci numbers).
%For $n=3$ on the Chebyshev nodes, $\D$ is (printing only at most five decimals)
%\begin{equation}
%    \D = \left[ \begin {array}{cccccccccccc}  0.0& 1.0& 0.0& 0.0& 0.0& 0.0&
% 0.0& 0.0& 0.0& 0.0& 0.0& 0.0\\ \noalign{\medskip} 0.0& 0.0& 2.0& 0.0&
% 0.0& 0.0& 0.0& 0.0& 0.0& 0.0& 0.0& 0.0\\ \noalign{\medskip}- 654.64&-
% 156.50&- 28.500& 704.0&- 144.0& 48.0&- 135.11&- 5.3333&- 16.0& 85.750
%&- 15.0& 1.5000\\ \noalign{\medskip} 0.0& 0.0& 0.0& 0.0& 1.0& 0.0& 0.0
%& 0.0& 0.0& 0.0& 0.0& 0.0\\ \noalign{\medskip} 0.0& 0.0& 0.0& 0.0& 0.0
%& 2.0& 0.0& 0.0& 0.0& 0.0& 0.0& 0.0\\ \noalign{\medskip}- 56.375&-
% 8.6250&- 0.75000& 45.111&- 26.0& 3.0& 26.0&{ 1.5187\times 10^{-13}}&
% 3.0&- 14.736& 2.5417&- 0.25000\\ \noalign{\medskip} 0.0& 0.0& 0.0&
% 0.0& 0.0& 0.0& 0.0& 1.0& 0.0& 0.0& 0.0& 0.0\\ \noalign{\medskip} 0.0&
% 0.0& 0.0& 0.0& 0.0& 0.0& 0.0& 0.0& 2.0& 0.0& 0.0& 0.0
%\\ \noalign{\medskip} 14.736& 2.5417& 0.25000&- 26.0&- 0.0&- 3.0&-
% 45.111&- 26.0&- 3.0& 56.375&- 8.6250& 0.75000\\ \noalign{\medskip}
% 0.0& 0.0& 0.0& 0.0& 0.0& 0.0& 0.0& 0.0& 0.0& 0.0& 1.0& 0.0
%\\ \noalign{\medskip} 0.0& 0.0& 0.0& 0.0& 0.0& 0.0& 0.0& 0.0& 0.0& 0.0
%& 0.0& 2.0\\ \noalign{\medskip}- 85.750&- 15.0&- 1.5000& 135.11&-
% 5.3333& 16.0&- 704.0&- 144.0&- 48.0& 654.64&- 156.50& 28.500
%\end {array} \right] \>.
%\end{equation}
In figure~\ref{fig:Dnorm} we find a log-log plot of the norms of $\D$ for these $n$. Remember that the degree of the interpolant is at most $3n+2$.  We see that the norm of $\D$ grows extremely rapidly for equally-spaced nodes (as we would expect).  For Chebyshev nodes there is still substantial growth (for confluency $3$; for confluency $2$ there is less growth, and for confluency $4$ there is more), but for $n=55$ and confluency $3$ at all nodes we have $\|\D\| $ approximately $10^{10}$ which still gives some accuracy in $\vec{Z}$.

\begin{figure}[h]
    \centering
    \includegraphics[width=0.5\textwidth]{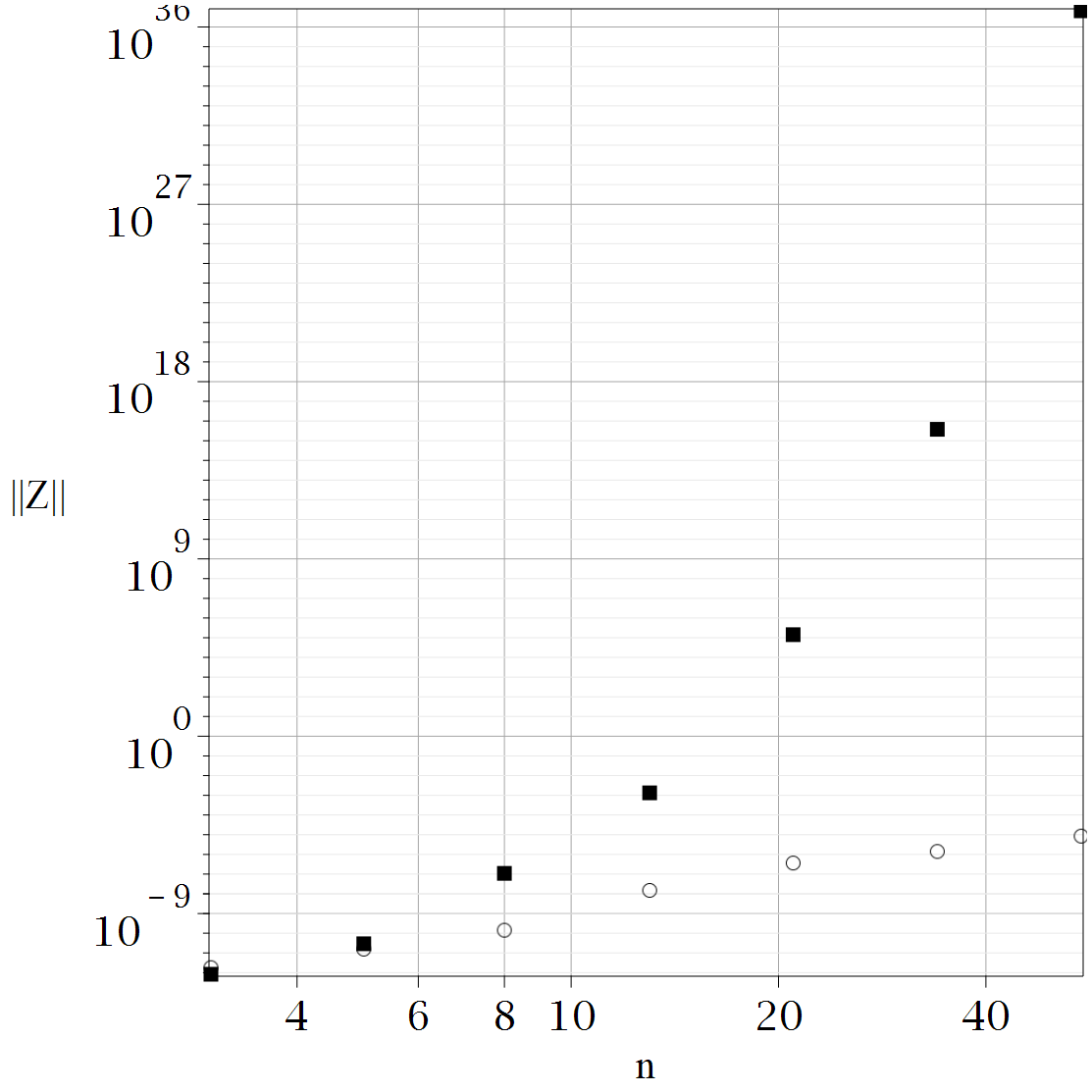}
    \caption{A comparison of norms of the vector $\vec{Z} = \D \mathbf{1}$ in Hermite interpolational bases. Equally-spaced nodes (solid box) and Chebyshev nodes (circles).  As expected, we have $\|\vec{Z}\| \approx \|\D\|\cdot 10^{-16}$ when working in double precision.}
    \label{fig:Znorm}
\end{figure}

In figure~\ref{fig:Znorm} we see the corresponding norms of $\vec{Z}$.  The behaviour is as predicted.

\par\medskip\noindent
\textbf{Remark}.
The confluency really matters.  If we use just simple Lagrange interpolation, that is confluency $s_i=1$ at each node, then the interpolation on $n=55$ Chebyshev nodes is in error by no more than $3.5\cdot 10^{-12}$.  Of course, the nominal degree is much lower than it was in the Hermite case with confluency $3$.  When we up the degree to $165$, the Lagrange error is no more than $1.5\cdot 10^{-11}$.  When the confluency is $3$, and $n=55$ which is comparable, the error is $1.4\cdot 10^{-5}$.

\section{Concluding Remarks}
Expressing a polynomial in a particular basis reflects a choice taken by a mathematical modeller.  We believe that choice should be respected. Indeed, changing bases can be ill-conditioned, often at least exponentially in the degree.  There are exceptions, of course: interpolation on roots of unity with a Lagrange basis can be changed to a monomial basis by using the DFT, and the conversion is perfectly well-conditioned; similarly changing from a Lagrange basis on Chebyshev-Lobatto points to the Chebyshev basis is also perfectly well-conditioned.  But, usually, one wants to continue to work in the basis chosen by the modeller. This is particularly true of the Bernstein basis, which has an optimal conditioning property: out of all bases that are nonnegative on the interval $[0,1]$, the Bernstein basis expression has the optimal condition number~\cite{Farouki(1996)}.  This property was extended to bases nonnegative on a set of discrete points by~\cite{Corless(2004)}, who proved that Lagrange bases can be better even than Bernstein bases.  See also~\cite{carnicer2017optimal}, who independently proved the same.

Differentiation is a fundamental operation, and it is helpful to be able to differentiate polynomials without changing bases.  This paper has examined the properties of the matrices for accomplishing this.  We found several of the results presented here to be surprising, notably that the Jordan Canonical Form for all the differentiation matrices considered here was the same. Likewise, that there is a uniform formula for a pseudo-inverse of all differentiation matrices of the type considered here was also a surprise.

One can extend this work in several ways.  One of the first might be to look at differentiation matrices for \textsl{compact finite differences}.  These are no longer always exact, and the matrices arising are no longer nilpotent (though they have null spaces corresponding to the polynomials of low enough degree that they are exact for). There are also some further experiments to run on the differentiation matrices we have studied in this paper already.  For instance, it would be interesting to know theoretically the growth of $\|\D^k\|$ for various dimensions $n$; we found that for the Bernstein basis of dimension $n+1$ we had $\| \D^n \|_\infty = 2^n n!$.  Since for the monomial basis we have $\| \D^n \|_\infty = n!$, this suggests that the natural scale for such a comparison is to divide by $n!$ and indeed that seems logical, because then in essence we are comparing the size of Taylor coefficients instead of comparing the size of derivatives and it is the Taylor coefficients that have a geometric interpretation in terms of location of nearby singularities.  We leave this study of the dependence of the norm of $\D^k$ in different bases for future work.

\section*{Acknowledgements}
This work was supported by a Summer Undergraduate NSERC Scholarship for the third author.  The second author was supported by an NSERC Discovery Grant.
We also thank ORCCA and the Rotman Institute of Philosophy.

\bibliographystyle{plain}
\bibliography{References,ref}

\end{document}